\documentstyle{amsppt} \magnification 1200
\input pictex
\UseAMSsymbols
\hsize 5.5 true in
\vsize 8.5 true in
\parskip=\medskipamount
\NoBlackBoxes

\def\mathbb{\Bbb}

\def\mathcal{\Cal}

\def\supp{\text{\rm supp\,}}

\def\Res{\text{\rm Res}}

\def\mod{\text{\rm {mod\,}}}
\def\ve{\varepsilon}
\def\vp{\varphi}

\def\mathbb{\Bbb}

\def\snint{\raise2pt\hbox{$_{^\not}$}\kern-3.5 pt\int}
\def\nint{\not\!\!\int}

\def\snint{\raise2pt\hbox{$_{^\not}$}\kern-3.5 pt\int}

\def\mes{\text {\rm mes\,}}
\def\vp{\varphi}
\TagsOnRight
\NoRunningHeads

\document
\parskip=\medskipamount
\topmatter
\title
On the correlation of the Moebius function with rank-one systems
\endtitle
\author
J.~Bourgain
\endauthor
\address
Institute for Advanced Study, Princeton, NJ 08540
\endaddress
\email
bourgain\@ias.edu
\endemail
\abstract
We explore the `Moebius disjointness property' in the special context of rank-one transformations and verify this phenomenon for many of the `classical'
models.
\endabstract
\endtopmatter

\noindent
{\bf (0). Introduction and Preliminaries}

This note is a follow up on \cite {B-S-Z}.
Recall that the general problem considered is the orthogonality of the Moebius sequence $\{\mu(n); n\in\Bbb Z_+\}$ with orbits of dynamical systems of zero
(topological) entropy.
To be more precise, as pointed out in \cite {B-S-Z}, we consider a uniquely ergodic topological model for the system.
Perhaps the simplest class of systems to study in this context are rank-one transformations and even in this class the problem seems highly nontrivial.
It was observed in \cite{B-S-Z} that if an orbit $f(n)=\vp(T^nx), n\in\Bbb Z_+$ is not orthogonal to the Moebius sequence, then there exist arbitrary
large $X\in\Bbb Z_+$ and subsets $ \Cal P_X\subset [p\in \Cal P; p\sim X]$ (where $\Cal P$ stands for the set of the primes) such that
$|\Cal P_X|\to \infty$ and for $p, q\in \Cal P_X$
$$
\overline{\lim_N}\, \Big| \frac 1N\sum^N_1 f(pn) \overline {f(qn)}\Big|>0\tag 0.1
$$
(this may be seen as a variant of the classical Vinogradov bilinear method to estimate sums involving $\mu(n)$ or the Van Mangoldt function $\Lambda(n)$).

It follows in particular that Moebius-orthogonality is implied by
 the disjointness of distinct powers $T^p$ and $T^q$ of the transformation $T$.
If $T$ is rank-one, then a result due to J.~King \cite{K} asserts that $T$ has the MSJ (minimal self-joining property) whenever $T$ is mixing, which in turn implies disjointness of
$T^p$ and $T^q$ for $p\not= q$ by a result of Del Junco and Rudolph \cite{D-R}.
Hence (0.1) certainly holds for mixing rank-one transformations, such as Ornstein's constructions and the Smorodinsky-Adams map.

Note also that there are examples of non-mixing rank-one transformations, for instance Chacon's transformation (see \cite{D-R-S}) that are MSJ.
On the other hand, there are natural classes of rank-one transformations (including generic interval exchanges, see \cite {V}) that are rigid and hence not MSJ.
Our goal here is to study the Moebius-orthogonality property directly, without relying on MSJ, for instance by establishing disjointness of $T^p$
and $T^q$, for certain $p\not= q$, by an ad hoc argument.
Our method is of spectral nature and is based on the observation that the spectral measure of a rank-one transformation is given by a generalized Riesz-product
(see \cite{B}). 
Formulating the disjointness of $T^p$ and $T^q$ from this perspective immediately leads to harmonic analysis problems related to the singularity of Riesz-products and we are able to treat
certain cases.
In order to make more precise statements, first recall the combinatorial definition of a rank-one transformation.

Following \cite {Fe}, a standard model of a rank one system is defined as follows.
We are given sequences of positive integers $w_n, n\in\Bbb Z_+$ and $a_{n, i}; n\in\Bbb Z_+, \break
1\leq i\leq w_n{-1}$ and define
$$
h_0=1, h_{n+1} = w_n h_n+\sum^{w_n-1}_{j=1} a_{n, j}.\tag 0.2
$$
Assume
$$
\sum^\infty_{n=1} \frac 1{w_nh_n} \Big(\sum^{w_n-1}_{j=1} a_{n, j} \Big) <\infty.\tag 0.3
$$
Define words $B_n$ on the alphabet $\{0, 1\}$ by
$$
B_0=0, \, B_{n+1}=B_n 1^{a_{n, 1}} B_n\cdots B_n 1^{a_{n, w_n-1}} B_n.\tag 0.4
$$
Consider the symbolic dynamical system $(X, T)$ where $X\subset \{0, 1\}^{\Bbb Z_+}$ consists of the sequence $(X_n)$ such that for every pair $s<t$,
$(x_s, \ldots, x_t)$ 
is a subsequence of some word $B_n$ and $T$ is the shift.

The `one' symbols in (0.4) between the words $B_n$ are called spacers.
If we assume
$$
\sup_{n, i} a_{n, i}< \infty\tag 0.5
$$
the topological system $(X, T)$ will be minimal (otherwise the only possible non-dense orbit is the sequence identically equal to one)
and uniquely ergodic.

\proclaim
{Theorem 1} Assume $T$ a rank-one system and (with the above notation)
$$
w_n< C, a_{n, j}<C.\tag 0.6
$$
Then $T$ satisfies the Moebius orthogonality property.
\endproclaim

As will be clear from the argument, the same conclusion may be reached under weaker assumptions then (0.6), but we certainly are not able to treat the
general case of a rank-one transformation at this point.

A few comments about the proof.
We start by assuming that $T$ is weakly mixing and later remove this hypothesis.
Assuming
$$
\min_{j\leq w_n} a_{n, j}< \max_{j\leq w_n} a_{n, j}\tag 0.7
$$
for sufficiently many values of $n$, we show that if $ \Cal P_0$ is a sufficiently large set of primes, we can find $p, q\in \Cal P_0$ 
such that $T^p$ and $T^q$ are disjoint, implying Moebius
orthogonality.

Otherwise (i.e. (0.7) only holds on a very thin set), the system $T$ satisfies a strong rigidity property in which situation another argument can be applied.
This is the same dichotomy as when estimating exponential sums
$$
\sum^N_1 \Lambda(n) e(n\theta) \text { or } \sum^N_{1} \mu(n) e(n\theta)\tag 0.8
$$
with $\big|\theta-\frac aq\big|<\frac 1{q^2}$, depending on whether $q$ is small or large (the first case relying on Dirichlet $L$-function theory and the second
on Vinogradov's bilinear method).

Next, we make a few observations and bring the generalized Riesz products into the picture.
Denote $\Bbb T=\Bbb R/\Bbb Z$ the circle, $e(\theta)= e^{2\pi i\theta}$.
Note that for $p, q$ distinct primes
$$
\int_0^1 \Big[\sum^N_1f(n) e(pn\theta)\Big] \overline{\Big[\sum^N_1 g(n)e(qn\theta)\Big]} d\theta= \sum_1^{N_1} f(qn) \overline{g(pn)} \text { with }
N_1 =\frac N{\max (p, q)}.\tag 0.9
$$
We introduce trigonometrical polynomials $(j\geq 1)$
$$
P_j(\theta) =\frac 1{\sqrt{w_j}}\sum_{k=0}^{w_j-1} e\big((kh_j+s_j(k))\theta\big)\tag 0.10
$$
$$
\text {where $s_j (0)=0, s_j(k)=a_{j, 1}+\cdots a_{j, k}$ for $1\leq k<w_j$}.
$$
Let $\vp\in\Cal C(X)$ be a continuous function on $X$, $\Bbb E[\vp]=0$.
By approximation, we may assume that $\vp$ depend on finitely many coordinates.

Following \cite{B}, if $x\in X$, $f(m)=\vp(T^m x)$ and $N=h_n$ large enough, we may write
$$
\frac 1{\sqrt N}\sum_1^N f(m) e(m\theta)= \omega(\theta)\prod^n_{n_0}  P_j(\theta)\tag 0.11
$$
with $\omega=\omega_\vp \in \Cal C(\Bbb T), \omega(\theta)=0$.

In general, given $\ve>0$ and $N$ large, the segment $I=[1, N]\cap \Bbb Z$ may be decomposed in segments of length $h_n\sim\ve N$ \big(we use here (0.2), (0.6)\big),
separated by a bounded number of spaces, with an initial and terminal segment of size $O(\ve N)$.
Hence, in view of (0.9), it will clearly suffice to prove that
$$
\int_{\Bbb T\backslash U}\prod^n_{n_0} | P_j(p\theta)| \ | P_j(q\theta)|d\theta\overset {n\to \infty} \to \longrightarrow 0\tag 0.12
$$
where $U\subset\Bbb T$ is a fixed neighborhood of 0.

Note that
$$
R_n(\theta)=\prod^n_1 | P_j(\theta)|^2\tag 0.13
$$
satisfies 
$$
\int_{\Bbb T} R_n(\theta) d\theta=1.\tag 0.14
$$
The maximal spectral type of $T$ is given by $\mu = $ weak$^*$ $\lim R_n(\theta) d\theta$ (cf. \cite {B}).

The proof of (0.12) occupies \S1.
As said earlier, we only establish this property under an additional assumption that excludes systems that are too rigid (see Lemma
2).
The remaining cases fall within a class of symbolic systems discussed in \S2.
For those, Moebius orthogonality is obtained in a different way, writing
$$
\Big|\sum^N_1\mu(n) f(n)\Big|\leq \int_0^1 \Big|\sum^N_1\mu(n) e(n\theta)\Big| \, \Big|\sum^N_1 f(n) e(n\theta)\Big| d\theta
\tag 0.15
$$
and obtaining good estimates for (restricted) $L^1$norms of the exponential sums $\sum^N_1f(n) e(n\theta)$ (see Theorem 7).
The class of systems introduced in \S2 turns out to capture also other classical rank one transformations such as generalized Chacon and Katok systems.

In fact, the bounds obtained in \S2 go beyond proving Moebius disjointness.
They may be combined with a standard Hardy-Littlewood type analysis involving minor and major arcs (a very crude version of this
technique already suffices for our purpose) in order to prove a pointwise prime number theorem for systems that fit the conditions
involved in \S2.
Natural examples of such systems are provided by three-interval exchanges transformations, satisfying certain (not to restrictive)
assumptions, as discussed in \S3. Thus the analysis from \S2 leads to the following result (cf. \cite{F-H-Z$_{1,2,3,4}$} for background).

\proclaim
{Theorem 8} Assume $T_{\alpha, \beta}$ a three-interval exchange transformation satisfying the Keane condition and such that the
associated three-interval expansion sequence $(n_k, m_k)_{k\geq 1}$ of integers fulfills the conditions
$$
\inf_k \frac {\min(n_k, m_k)}{n_k+m_k}>0
$$
and
$$
\min(n_k, m_k)>C_0
$$
with $C_0$ a sufficiently large constant.

Then $T_{\alpha, \beta}$ has the Moebius disjointness property and also satisfies a (pointwise) prime number theorem.
\endproclaim

The harmonic analysis involved in \S1 and \S2 may be of independent interest and deserving to be further explored.

\bigskip

\noindent
{\bf (1). Estimates on Riesz-Products}

In this section, we develop a method to estimate the left side of (0.12).
For simplicity, set $n_0=1$. Let
$$
W=\max_n w_n \text { and } A=\max_{n, i} a_{n, i}<\infty.
$$
Clearly, from (0.2)
$$
|h_{j+1} - w_jh_j|< A w_j
$$
and iterating
$$
\align
|h_{j+2} - w_{j+1} w_j h_j&|< A w_{j+1}+Aw_{j+1} w_j\\
&\vdots
\endalign
$$
For $n>j$
$$
\align
|h_n-w_{n-1}\cdots w_j h_j|&< A w_{n-1}+Aw_{n-1} w_{n-2}+\cdots +Aw_{n-1}\ldots w_j\\
&<2A w_j\cdots w_{n-1}.\tag 1.1
\endalign
$$
Let $n_0<n$ and write
$$
\prod_{j\leq n}  P_j =\prod_{j< n_0}  P_j . \prod_{n_0\leq j\leq n} P_j = F.G
$$
where, by (0.2)
$$
\align
\supp \hat F&\subset \sum_{j<n_0}\{ k h_j+ s_j (k); 0\leq k\leq w_j - 1\}\\
&\subset [\theta, h_{n_0}]\tag 1.2
\endalign
$$
and by (1.1),
$$
\spreadlines{6pt}
\align
\supp \hat G &\subset\sum_{n_0\leq j\leq n}\{kh_j+s_j(k), 0\leq k\leq w_j-1\}\subset\\
&\{ k h_{n_0}; 0\leq k\leq w_{n_0} -1\}\\
&+\sum_{n_0< j\leq n}\{kw_{j-1} \cdots  w_{n_0} h_{n_0}; 0\leq k\leq w_j-1\}\\
&+[0, 2WAw_{n_0}\cdots w_n+(n-n_0)AW].\tag 1.3
\endalign
$$
Thus $t\in\supp\hat G$ is of the form
$$
t=t_0+(t-t_0)
$$
with
$$
t_0\in\Lambda=\Lambda_{n_0, n} =\sum_{n_0<j\leq n} \{ kw_{j-1} \ldots w_{n_0} h_{n_0}; 0\leq k< w_j\}\tag 1.4
$$
and
$$
|t-t_0|< 3AW^{n-n_0+2}.\tag 1.5
$$
Assume
$$
h_{n_0}\gg W^{n-n_0} \text { (which amounts to $n-n_0< cn$ for some constant $c>0 $)}\tag 1.6
$$

Let $0<\delta \ll1$ take $\psi\in\Bbb T$ satisfying
$$
\Vert\psi\Vert< \delta h^{-1}_{n_0}.\tag 1.7
$$
Denote for $n_0\leq j\leq n$ by
$$
P_j(\theta, \psi)=\frac 1{\sqrt {w_j}} \sum_{k=0}^{w_j-1} e\big(( k h_j+s_j(k))\theta\big) e(kw_{j-1} \cdots w_{n_0} h_{n_0}\psi).\tag 1.8
$$
It follows from the preceding that if 
$$
\xi=\xi'+t\in \supp \hat F+\supp \hat G
$$
and $\psi$ satisfies (1.7), then
$$
\align
|e(\xi\psi)- e (t_0\psi)|&\leq |e(\xi'\psi)-1|+ |e\big((t-t_0)\psi\big)-1|\\
&\underset{(1.2), (1.5)}\to\leq h_{n_0} \Vert\psi\Vert+C W^{n-n_0} \Vert\psi\Vert\\
&\underset {(1.1), (1.7)}\to < \ 2 \delta.\tag 1.9
\endalign
$$
Therefore clearly
$$
\Big\Vert \prod_1^n  P_j(\theta+\psi) -\prod_{j<n_0}  P_j (\theta) \prod_{n_0\leq j\leq n}  P_j (\theta,\psi)\Big\Vert_{{L^2}(d\theta)} 
\lesssim \delta\tag 1.10
$$
and hence, by Cauchy-Schwartz inequality and (0.14)
$$
\Big\Vert R_n (p\theta + p\psi)-R_{n_0-1} (p \theta)\prod_{n_0\leq j\leq n} |P_j(p\theta, p\psi)|^2\Big\Vert_{L^1(d\theta)} \leq \delta p.\tag 1.11
$$
Assume the left side of (0.12) larger than $c_0$.
Taking $\delta\sim \frac {c_0^2}{p+q}$, it follows from (1.11) that (denoting $\nint$ the normalized integral)
$$
\frac {c_0}2< \nint_{\Vert\psi\Vert< \delta h^{-1}_{n_0}} \int_{\Bbb T\backslash U} R_{n_0-1} (p\theta)^{1/2} 
R_{n_0-1} (q\theta)^{1/2}  \prod_{n_0\leq j\leq n} |P_j(p\theta, p\psi)| \, |P_j
(q\theta, q\psi)| d\theta d\psi.\tag 1.12
$$
Choose $n_1>n_0$ such that $w_{n_1-1}\cdots w_{n_0} >\frac 1\delta$, thus
$$
n_1 -n_0\lesssim \log\frac {p+q}{c_o}.\tag 1.13
$$
Since $|P_j|\leq \sqrt{w_j}$, it follows from (1.12) that
$$
\align
\nint_{\Vert\psi\Vert<\delta h^{-1}_{n_0}} \int_{\Bbb T\backslash U}  R_{n_0-1} (p\theta)^{1/2} 
R_{n_0-1}(q\theta)^{\frac 12}
 \prod_{n_1\leq j\leq n} |P_j(p\theta, p\psi)|\, |P_j(q\theta, q\psi)| d\theta d\psi\\
>\frac 12 \frac \delta W c_0 =c_1\sim c_0^3\tag 1.14
\endalign
$$
(the range of the primes $p, q\cdots$ is considered bounded).

Next, making a change of variable $\psi' =w_{n_1-1}\cdots w_{n_0} h_{n_0}\psi$ and defining for $n_1\leq j\leq n$
$$
 P_j(\theta, \psi')= \frac 1{\sqrt{w_j}} \sum_{k=0}^{w_j-1} e\big((kh_j+s_j(k))\theta\big) e(kw_{j-1} \cdots w_{n_1}\psi')\tag 1.15
$$
we obtain
$$
\int_{\Bbb T\backslash U} \int_{\Bbb T} R_{n_0-1} (p\theta)^\frac 12 R_{n_0-1}(q\theta)^{\frac 12}
\prod_{n_1\leq j\leq n}|P_j(p\theta , p\psi')| \ 
|P_j(q\theta, q\psi')| d\theta d\psi'>c_1.\tag 1.16
$$
We will exploit the $\psi'$-variable.
By (1.15), for $j=n_1$
$$
 P_{n_1} (\theta, \psi')=\frac 1{\sqrt{w_{n_1}}} \sum^{w_{n_1}-1}_{k=0} e \big((kh_{n_1}+s_{n_1}(k)\big)\theta\big)
e (k\psi').\tag 1.17
$$
Making a shift $\psi' \mapsto\psi'+\frac r{w_{n_1}} (0\leq r< w_{n_1})$, it follows from (1.15) that the  factors
$ P_j(\theta, \psi'), j\geq n_1$, are preserved.

Define
$$
\rho_{n_1}(\theta) =\max_{\psi'} \frac 1{w_{n_1}} \sum_{r=0}^{w_{n_1}-1} \Big| P_{n_1} \Big(p\theta, p\psi'+\frac {pr}{w_{n_1}}\Big)
\Big| \, \Big|  P_{n_1} \Big(q\theta, q\psi' + \frac {qr}{w_{n_1}}\Big) \Big|.\tag 1.18
$$
Then
$$
(1.18) \leq\int_{\Bbb T\backslash U} \int_{\Bbb T} R_{n_0-1} (p\theta)^{\frac 12} R_{n_0-1}(q\theta)^{\frac 12} \rho_{n_1}(\theta) 
\prod_{n_1< j\leq n} | P_j(p\theta, p\psi')| \, | P_j(q\theta, q\psi')|d\theta d\psi'.\tag 1.19
$$
Defining in general
$$
\rho_j(\theta)= \max_{\psi'} \frac 1{w_j} \sum_{r=0}^{w_j-1} \Big| P_j\Big(p\theta, p\psi'+\frac {pr}{w_j}\Big)\Big| \, 
\Big| P_j\Big(q\theta, q\psi' + \frac {qr}{w_j}\Big)\Big|.\tag 1.20
$$
Iteration implies that
$$
C_1\leq \int_{\Bbb T\backslash U}  R_{n_0-1} (p\theta)^{\frac 12} R_{n_0-1} (q\theta)^{\frac 12} \prod^n_{j=n_1} \rho_j(\theta) d\theta.\tag 1.21
$$
Assume
$$
p, q> W.\tag 1.22
$$
We analyze $\rho_j$.
Fix $\psi'$ and write (setting $v=w_j, P= P_j)$
$$
(1.20) =\frac 1{2v} \sum^{v-1}_{r=0} \Big| P\Big(p\theta, p\psi'+\frac {pr}v\Big)\Big|^2
+ \frac 1{2v} \sum^{v-1}_{r=0} \Big | P\Big (q \theta, q{\psi'}+\frac {qr} v\Big)\Big|^2
$$
$$
\align
&-\frac 1{2v} \sum^{v-1}_{r=0}\Big| \ \Big| P\Big(p\theta, p\psi' +\frac {pr}v\Big)\Big|-
\Big | P\Big(q \theta, q{\psi'}+ \frac {qr}v\Big)\Big|\ \Big|^2\\
&\overset {(1.21)}\to = 1 -\frac 1{2v}\sum^{v-1}_{r=0} \Big|\ \Big| P\Big(p\theta, p\psi'+\frac {pr}{v}\Big)\Big| 
-\Big| P\Big(q\theta, q\psi'+\frac {qr}v\Big)\Big|^2\Big|.\tag 1.23
\endalign
$$
Also
$$
\align
&\frac 1v \sum^{v-1}_{r=0} \Big| \ \Big| P\Big(p\theta, p\psi'+\frac {pr}v\Big)\Big|^2 -
\Big| P\Big(q\theta, q\psi'+ \frac {qr}v\Big) \Big|^2\Big|\leq\\
&\Big\{\frac 1v \sum^{v-1}_{r=0} \big| \,|\ | - |\ |\,\big|^2 \Big\}^{\frac 12}. \Big\{\frac 1v \sum^{v-1}_{r=0} 
\big[ | \ |+| \ |\big]^2\Big\}^{\frac 12}\leq\\
&2\Big\{\frac 1v \sum^{v-1}_{r=0} \big| \, |\ | - |\ | \,\big|^2\Big\}^{\frac 12}
\endalign
$$
and thus
$$
\align
(1.23) &< 1-\frac 1{8v^2} \Big\{ \sum^{v-1}_{r=0} \Big| \ \Big| P \Big(p\theta, p\psi'+\frac {pr}v\Big)\Big|^2 -\Big | \
\Big| P\Big( q\theta, q\psi'+\frac {qr}v\Big)\Big|^2 \Big| \Big\}^2\\
& \leq 1-\frac 18 (1.25)^2\tag 1.24
\endalign
$$
where by (1.15) and change  of variables $\eta =w_{j-1} \cdots w_{n_1}\psi'+ h_j\theta$,
$$
(1.25) =\min_\eta \frac 1{v^2} \sum^{v-1}_{r=0} \Big|\ \Big|\sum^{v-1}_{k=0} e\Big(kp\eta+s_j (k)p\theta+ k\frac {pr}v\Big)\Big|^2 -
\Big|\sum^{v-1}_{k=0} e\Big(kq\eta+s_j(k) q\theta +k\frac {qr}v\Big)\Big|^2\Big|.\tag 1.26
$$
Assume moreover
$$
p\equiv 1\equiv q (\mod v).\tag 1.27
$$

Note that since $v\in \{w_n\}$ which  is a finite set, we may always ensure (1.27) for all $v$ by restriction of the 
primes from the set $\Cal P$ in (0.1).
\medskip

Fix $\eta$ in (1.26) and minorize by
$$
\frac 1{v^3}\Big|\sum^{v-1}_{r=0} e\Big(\frac {-r}v\Big) \Big\{\Big|\sum^{v-1}_{k=0} e(k p\eta+s_j(k) p\theta+\frac {kr}v\Big)\Big|^2 -
\Big|\sum^{v-1}_{k=0} e(kq\eta+s_j(k) q\theta +\frac {kr}v\Big)\Big|^2\Big\}\Big|
$$
which by (0.11) equals (setting $a(k)=a_j (k), a(0)=0$)
$$
\align
&\frac 1{v^3}\Big|e(p\eta) \Big(\sum^{v-1}_{k=1} e(a(k) p\theta)\Big) + e\big(-(v-1)p\eta\big) e\Big(-\Big(\sum^{v-1}_{1} a(k)\Big)p\theta\Big)\\
&- e(q\eta) \Big(\sum^{v-1}_{1} e(a(k) q\theta)\Big) 
-e\big(-(v-1) q\eta\big) e\Big(-\Big(\sum^{v-1}_{1} a(k)\Big) q\theta\Big)\Big|\tag 1.28
\endalign
$$
the crucial fact in our analysis is the following property

\noindent
(1.29). 

For all $\ve>0$, there is $\ve_1> 0$ such that
$$
\mes[\theta; \min_\eta (1.28) <\ve_1]<\ve.
$$
If (1.29) fails, then for some $\ve>0$
$$
\mes[\theta; \min_\eta (1.28) =0] >\ve.\tag 1.30
$$
Introduce the rational function $f(x, y) \in\Bbb C(X, Y)$ defined as
$$
f(x, y) =y^p\sum^{v-1}_{1} x^{a(k)p} +y^{-(v-1)p} x^{-p\sum a(k)} -y^q \sum_1^{v-1} x^{a(k)q} -y^{-(v-1)q} x^{-q\sum a(k)}.\tag 1.31
$$
It follows from (1.30) that for $x\in [e^{i\theta}; \theta \in\Omega]$, $\mes\Omega>0$, the equations
$$
 \qquad\qquad\qquad
\Bigg\{
\aligned
&f(x, y)=0\qquad\qquad\qquad\qquad \  \qquad\qquad\qquad\qquad\qquad \ (1.32) \\ 
&f\Big(\frac 1x, \frac 1y\Big) =0 \qquad\qquad\qquad\qquad\qquad\qquad\qquad\qquad\quad \  \ (1.33)
\endaligned
$$
have a common root in $y$.

Assume $p<q$. Then
$$
g_1 (x, y) =y^{(v-1)q} f(x, y) \text { and } g_2(x, y)=y^q f\Big(\frac 1x, \frac 1y\Big)
$$
are polynomials in $y$ and, from the preceding, their resultant
$$
\Res_y \big(g_1(x, y), g_2 (x, y)\big)\in\Bbb C(X)
$$
vanishes identically.
Therefore (1.32), (1.33) have a common component $\Gamma$.

We show that this if impossible, provided
$$
a_+=\max \big(a(1), \ldots , a(v-1)\big) > a_- =\min \big(a(1), \ldots, a(v-1)\big).\tag 1.34
$$
Using $x$ as a local coordinate on $\Gamma$, it follows from (1.31) that
$y\to \infty$ for $x\to 0$ and hence we obtain Puisseux expansion
$$
y=\xi x^{-\alpha}+\sum_{\beta>-\alpha} c_\beta x^\beta \qquad (\xi\not= 0, \alpha\in \Bbb Q_+)\tag 1.35
$$
solving both (1.32), (1.33).

Since $p<q$ and $\alpha>0$, we derive from (1.31), (1.34), (1.35) (distinguishing the cases $\alpha>a_-$, $\alpha<a_-, \alpha =a_-)$ that
$$
a_--\alpha =\alpha (v-1)-\sum^{v-1}_{1} a(k)
$$
hence
$$
v\alpha = \sum_1^{v-1} a(k)+a_- .\tag 1.36
$$
Similarly, from (1.33), we obtain that
$$
v\alpha =\sum_1^{v-1} a(k)+a_+.\tag 1.37
$$
Hence $ a_+=a_-$, contradicting  assumption (1.34).

This proves (1.29). 
\bigskip

Returning to (1.24), we proved that if
$$
\max_{0<k< w_j} a_j(k)> \min_{0<k<w_j} a_j(k)\tag 1.38
$$
then
$$
\rho_j(\theta) < 1 - \frac 1 8 \ve_1^2 \ \text { where } \ \ve_1=\ve_1(\ve)>0\tag 1.39
$$
provided
$$
\theta\not\in \Omega_j\tag 1.40
$$
where
$$
\Omega_j =[\theta\in\Bbb T; \min_\eta (1.25) < \ve_1] \text { satisfies } \mes \Omega_j <\ve\tag 1.41
$$
with $v=w_j$ and $a(k) =a_j(k)$ in (1.25).

Here $j\in [n_1, n]$.
Recall that by (1.6), (1.13) we may take $n-n_1\sim n$.
Assume (1.38) for all $j$ (or at least sufficiently many).
We may then specify $v$ and a configuration $\big(a(k)\big)_{0<k<v}$ such that for all $j\in\frak S$, where
$$
|\frak S|> (n-n_1)\frac 1{WA^{W-1}}\sim n\tag 1.42
$$
we have $w_j=v$ and $a_j(k) =a(k)$ for $j\in\frak S$.
Thus $\rho_j(\theta)=\rho(\theta)$ for $j\in\frak S$.

Next, fix $\ve$ and take $\ve_1$ small enough to ensure
$$
\mes\Omega<\ve.\tag 1.43
$$

Since $\rho_j(\theta)\leq 1$, we obtain 
$$
\align
c_1 \leq (1.21) &\leq \int_{\Bbb T\backslash U} R_{n_0-1} (p\theta)^{1/2}
R_{n_0-1}(q\theta)^{1/2} \prod_{j\in\frak S} \rho_j(\theta)d\theta\\
&\overset{(1.39)} \to < \big(1-\frac{\ve_1^2}8\Big)^{|\frak S|} +\int_{\Omega\backslash U} R_{n_0-1}(p\theta)^{1/2} R_{n_0-1}(q\theta)^{1/2}d\theta.
\endalign
$$
Letting $n\to \infty$ along an appropriate subsequence, we find
$$
\mu_1(\bar\Omega\backslash U) \geq \lim_n \int_{\Omega\backslash U} R_n(p\theta)^{1/2} R_n(q\theta)^{1/2} d\theta\geq c_1\tag 1.44
$$
with $\mu_1$ a weak$^*$-limit point of $\{R_n(p\theta)^{1/2} R_n(q\theta)^{1/2}\} \subset L^1 (\Bbb T)$.

We assume that $T$ is weakly mixing, hence the spectral measure $\mu$ has no atoms outside $0$.
Since $p, q$ are distinct primes, it follows that also $\mu_1$ has no atoms outside $0$.

Note that (1.28) =$|Q(\theta, \eta)|$, where $Q(\theta, \eta)$ is a fixed trigonometric polynomial on $\Bbb T^2$.
Introduce real variables
$$
\Bigg\{\aligned t_1=cos \, \theta\\ t_2= \sin\theta\endaligned
\qquad \Bigg\{\aligned u_1=\cos\eta\\ u_2= \sin \eta\endaligned 
$$
and write $|Q(\theta, \eta)|^2 =Q_1(t_1, t_2, u_1, u_2)\in\Bbb R[t_1, t_2, u_1, u_2]$.

The set $\Omega$ is then obtained as a projection Proj$_{t_i} (\Omega')$ with
$$
\Omega'= \{(t_1, t_2, u_1, u_2)\in [0, 1]^4; t^2_1+t_2^2 =1, u_1^2+u_2^2 =1 \text { and } Q_1 (t_1, t_2, u_1, u_2)<\ve^2_1\}.\tag 1.45
$$
From semi-algebraic set theory, it follows that $\Omega$ consists of at most $C_1=C_1(A, W, p, q)$ intervals in $\Bbb T$
of total measure at most $\ve$, by (1.43).

But, letting $\ve\to 0$ (1.44) implies that $\mu_1$ has atoms in $\Bbb T\backslash U$, a contradiction.

Returning to condition (1.38) and previous argument, it clearly suffices in fact to have (1.38) satisfied for 
$j\in\frak S_1 \subset [n_1, n]$ where $|\frak S_1|> C(\ve, p, q, A, W)=C$.
Assume that this is not the case, i.e. there is $\frak S_1\subset[n_1, n] , |\frak S_1|<C$ such that
$$
a_j(1)=\cdots =a_j(w_j-1) \text { for } 0\leq k<w_j \text { if } j\in [n_1, n]\backslash \frak S_1.\tag 1.46
$$
Thus for $j\in [n_1, n] \backslash \frak S_1$, the word $B_{j+1}$ has the form
$$
B_{j+1} =\underbrace{B_ja_j(1)\cdots B_ja_j(1)}_{w_j-1} B_j.\tag 1.47
$$
Considering $B_{j+2}$ for $j,j+1 \in [n_1, n]\backslash\frak S_1$, we obtain
$$
\align
B_{j+2}& = B_{j+1} a_{j+1}(1) B_{j+1} \cdots  a_{j+1}(1) B_{j+1}\\
&=B_j a_j(1) B_ja_j(1)\cdots B_j a_{j+1}(1) B_ja_j(1) B_ja_j(1) \cdots  B_j a_j(1)B_j\tag 1.48
\endalign
$$
and the spacer condition (1.38) will hold unless
$$
a_j(1)=a_{j+1}(1).\tag 1.49
$$
Again assume (1.49) for all $j\in [n_1, n]\backslash \frak S_2$ with $\frak S_1\subset\frak S_2, |\frak S_2|< C$.

Clearly the set $[n_1, n] \backslash \frak S_2$ is the union of at most $C$ intervals $J_\alpha = ] m_\alpha, n_\alpha]$, for which, by (1.49), the
word $B_{n_\alpha+1}$ has the form
$$
B_{n_\alpha +1}= \underbrace {B_{m_\alpha} a_{m_\alpha}(1) B_{m_\alpha} a_{m_\alpha}(1) \cdots B_{m_\alpha} a_{m_\alpha}(1)
B_{m_\alpha}} _{w_{m_\alpha} w_{m_\alpha+1} \cdots w_{n_\alpha}}. \tag 1.50
$$
Summarizing, what we proved is the following.

Suppose
$$
\int_{\Bbb T\backslash U} \prod^n_1 | P_j(p\theta)| \, | P_j(q\theta)| d\theta>c>0\tag 1.51
$$
(we use the notation $c, C$ for various quantities independent on $n$).

Then there is $n_0<n, n-n_0>cn$ where the word $B_n$ is obtained from $B_{n_0}$ by a system $W_0=B_{n_0}, W_1, \ldots, W_r=B_n$ where $W'=W_s$
relates to $W=W_{s-1}$ by a formula
$$
W' =W^k 1^{b_1} W^k 1^{b_2}\ldots W^k 1^{b_\ell}\tag 1.52
$$
for some $k, \ell, b_1, \ldots, b_\ell$ and $r, \ell, b_i<C$.

Next, we aim to make a similar statement with $n$ replaced by $n_0$.

Proceed as follows.
Let $n_0$ be arbitrary (large) and $n>n_0$.
Iterating (0.4), the word $B_n$ may be written in the form
$$
B_n=B_{n_0} 1^{a_1} B_{n_0} 1^{a_2}\cdots 1^{a_\ell} B_{n_0}\tag 1.53
$$
with $\ell =w_{n_0}\ldots w_n -1 $ and $a_1, \ldots, a_\ell\leq A$.

Hence
$$
\align
\prod^n_1 P_j(\theta) &= \Big[\prod^{n_0}_1  P_j(\theta)\Big] \Big[ \frac 1{\sqrt{\ell+1}} \Big(1+\sum^\ell_{m=1} e\big((mh_{n_0}+
a_1+\cdots+a_m)\theta\big)\Big)\Big]\\
&=\Big[\prod_1^{n_0}  P_j(\theta)\Big]. Q(\theta).\tag 1.54
\endalign
$$
Introduce again an additional variable $\psi$,
$$
|\psi| <\delta h_{n_0}^{-1}.\tag 1.55
$$
For fixed $\psi$ satisfying (1.55)
$$
\Big\Vert \prod^n_{1}  P_j (\theta+\psi)-\prod^{n_0}_{1}  P_j(\theta). Q(\theta+\psi)\Big\Vert_{L^2(d\theta)}\lesssim \sqrt\delta
$$
and
$$
\Big\Vert\prod^n_{1}(p\theta+p\psi)-\prod^{n_0}_{1}  P_j(p\theta)Q(p\theta+p\psi)\Big\Vert_{L^2(d\theta)} \lesssim \sqrt{p\delta}.\tag 1.56
$$
Hence, taking $\delta\sim\frac{c_0^2}{p+q}$, we obtain
$$
\align
&c_0 <\int\prod^n_1| P_j(p\theta)| \, | P_j(q\theta)|d\theta =\nint_{|\psi|<\delta h^{-1}_{n_0}} \int\prod^n_{1}
| P_j(p\theta+p\psi)| \, | P_j (q\theta+ q\psi)|d\theta d\psi\\
&< \nint_{|\psi| <\delta h_{n_0}^{-1}} \int\prod^{n_0}_1 | P_j(p\theta)| \, | P_j(q\theta)| \ |Q(p\theta+p\psi)| \, 
|Q(q\theta+q\psi)|d\theta d\psi+\frac {c_0}2.\tag 1.57
\endalign
$$
Next, the $h_{n_0}$-separation of the $Q$-frequencies allows to estimate
$$
\nint_{|\psi|<\delta h^{-1}_{n_0}} \Big|Q(p\theta+ p\psi)\Big|^2 d\psi<\frac 1\delta\int^1_0 \Big|Q\Big(p\theta+ p\frac {\psi'}{h_{n_0}}\Big)
\Big|^2 d\psi'\lesssim\frac 1\delta.\tag 1.58
$$
From (1.57), (1.58) and Cauchy-Schwarz, it follows that
$$
\int\prod^{n_0}_1 | P_j(p\theta)| \, | P_j(q\theta)|d\theta\gtrsim c_0 \delta\sim c_0^3.\tag 1.59
$$
We may then repeat (1.52) with $n$ replaced by $n_0$.

Iteration finally leads to the following statement

\proclaim
{Lemma 2} Let $p\not= q$ be fixed primes such that $p, q>W$ and $p\equiv 1\equiv q$ $(\mod w_j)$ for all $j$. Assume
$$
\int\prod^n_1 | P(p\theta)| \, | P_j(q\theta)|d\theta>c\tag 1.60
$$
for some constant $c>0$.
Then $B_n$ may be obtained from a symbolic system $W_0, W_1, \ldots, W_r =B_n$ with $|W_0|<C$, $r< C\log n$ and $W'=W_s$  obtained from
$W=W_{s-1}$ as
$$
W'= W^k 1^{b_1} W^k 1^{b_2}\cdots W^k 1^{b_\ell}\tag 1.61
$$
with $\ell, b_i<C$.
Moreover $|W_{s-C}|=o(|W_s|)$.
\endproclaim

In this situation of `strong rigidity', the Moebius orthogonality will be established using a different argument, $n.\ell.$ expressing
$$
\sum^N_{1} \mu(n) f (n) =\int\Big[ \sum^N_1\mu(n) e(n\theta)\Big] \overline {\Big[\sum^N_1 f(n) e(n\theta)\Big]} d\theta\tag 1.62
$$
and bounding the exponential sums.
This will be pursued in the next section, in a greater generality (of independent interest).

Finally, let us also explain how to remove the assumption that $T$ is weakly mixing.
Recall the classical estimate
$$
\Big\Vert\sum^N_1\mu(n) e(n\theta)\Big\Vert_{L^\infty (\Bbb T)} \underset A\to \ll N(\log N)^{-A}.\tag 1.63
$$
This allows us to assume that the function $\vp$ in (0.11) satisfies
$$
\frac 1M\sum_1^M\vp (T^mx) e(m\theta)=o(1)\tag 1.64
$$
for $\theta \in\Cal E$ (=an arbitrary given finite set) and $M$ large enough (we subtract the projection on possible eigenstates of the
system $(X, T)$).
Therefore (by a further perturbation) the function $\omega$ in (0.11) may be taken to satisfy $\omega(\theta)=0$ for $\theta \in\Cal E$ and in
(0.12), $U$ may be taken to be a neighborhood of $\Cal E$.

In particular, suitable choice of $U$ will make (1.44) impossible for $\Omega$ the union of a bounded number of $\ve$-intervals with $\ve\to 0$.

\bigskip

\noindent
{\bf (2). A Class of Symbolic Systems}

Consider a symbolic system on the alphabet $0, 1$ with order-$n$ words $W\in \Cal W_n$ of the form
$$
W=W_1^{k_1} W_2^{k_2} \cdots W_r^{k_r} \text { for some } W_1, \ldots, W_r\in \Cal W_{n-1}' =\bigcup_{m<n} \Cal W_m\tag 2.1
$$
where we assume that $r$ remains uniformly bounded $r<C$. 

We also assume the following property for the system $\{\Cal W_n\}$.

Let $W\in\Cal W_n$ and express $W$ in words $W'\in\Cal W_{n-s}'$, $0<s\leq n$, by iteration of (2.1). Then
$$
\frac {|W|}{\max |W'|}>\beta(s)\tag 2.2
$$
where
$$
\frac {\log\beta(s)}{s} \overset{s\to\infty}\to \longrightarrow \infty.\tag 2.3
$$

\noindent
{\bf Remark.}
In fact, in the sequel we will only use the property that \hfill\break $\beta(s)>C_0 s$, $s$ large, for some sufficiently large constant $C_0$.

To a given word $W=(x_1, \ldots, x_\ell)$, we associate the trigonometric polynomial
$$
 P_W(\theta) =\sum^\ell_{1} x_m e(m\theta).\tag 2.4
$$
Hence, if $|W_i|=\ell_i$ in (2.1)
$$
\align
 P_W(\theta)=& P_{W_1}(\theta) \Big[\sum^{k_1-1}_{j=0} e(j\ell_1\theta)\Big]+\\
& P_{W_2} (\theta) \Big[\sum_{j=0}^{k_1-1} e(j\ell_2 \theta)\Big] e(k_1 \ell_1 \theta)+\cdots +\\
& P _{W_r}(\theta)\Big[ \sum_{j=0}^{k_r-1} e(j\ell_r\theta)\Big] e\big((k_1\ell_1+\cdots+k_{r-1} \ell_{r-1})\theta\big).\tag 2.5
\endalign
$$

In order to bound $L^1$-norms, we rely on the following

\proclaim
{Lemma 3} Let $I\subset\Bbb T$ be an interval of size $\delta$ and $\ell\gtrsim \frac 1\delta$, $|W|=\ell$. Then
$$
\align
&\int_I | P_W(\theta)|\Big|\sum^k_{j=0} e(j\ell\theta)\Big| d\theta<\\
&C\log (k+2).\int_{I'} | P_W(\theta)|d\theta+ C\log(k+2) e^{-c(\log \ell)^{4/3}}\tag 2.6
\endalign
$$
where $I'$ is a $\frac {(\log \ell)^2}\ell$-neighborhood of $I$.
\endproclaim
\medskip

\noindent
{\bf Proof.} Estimate
$$
\align
\int_I| P_W(\theta)|\Big|\sum^k_{j=0} e(j\ell\theta) \Big|d\theta &\lesssim \int_I | P_W(\theta)| 
\Big[\frac 1{\Vert\frac{\ell\theta} 2\Vert +\frac 1k}\Big]d\theta\\
&\lesssim \int_I d\theta \Big[\max_{|\psi|<\frac 1 \ell} | P_W(\theta+\psi)| \Big] \ \Big[ \int\frac {d\theta}{\Vert\frac \ell 2
\theta +\psi\Vert+\frac 1k}\Big]\\
&\lesssim \log (2+k)\int_I \max_{|\psi|<\frac 1\ell} | P_W (\theta+\psi)|d\theta.\tag 2.7
\endalign
$$
Let $K$ be a trigonometric polynomial satisfying
$$
\hat K = 1\text { on } [-\ell, \ell]\tag 2.8
$$
$$
\supp\hat K\subset [-2\ell, 2\ell]\tag 2.9
$$
$$
|K(\theta)|\lesssim \ell\exp \big(-c(\ell\Vert\theta\Vert)^{2/3}).\tag 2.10
$$
By (2.8), $ P_W= P_W*K$.
Hence
$$
\align
\max_{|\psi|<\frac 1\ell} | P_W(\theta +\psi)|&\leq\int | P_W(\eta)|\Big[\max_{|\psi|<\frac 1\ell} |K(\theta+\psi-\eta)|\Big]d\eta\\
&\overset{(2.10)}\to <\int_{\Vert \theta-\eta\Vert<\frac {(\log \ell)^2}{\ell}} | P_W(\eta)|
\big[\max_{|\psi|<\frac 1\ell}|K(\theta+\psi-\eta)|\big]d\eta\\
&+\ell e^{-c(\log\ell)^{4/3}}.\tag 2.11
\endalign
$$
Also, from (2.10)
$$
\align
\int\Big[\max_{|\psi|<\frac 1\ell} |K(\theta+\psi-\eta)|\Big]d\theta& \lesssim \ell \int\max_{|\psi|<\frac 1\ell} e^{-c(\ell\Vert \theta+
\psi-\eta\Vert)^{2/3}} d\theta\\
&\lesssim \ell \int e^{-c(\ell\Vert\theta\Vert)^{2/3}} d\theta =O(1)
\endalign
$$
and from (2.11)
$$
\int_I\max_{|\psi|<\frac 1\ell}|P_W(\theta+\psi)|d\theta\lesssim \underset{I+[-\frac {(\log \ell)^2}\ell, \frac
{(\log \ell)^2}\ell ]}\to {\int|P_W|} +e^{-c(\log \ell)^{4/3}}\tag 2.12
$$
proving (2.6).

Taking $I=\Bbb T$, (2.6) becomes
$$
\int| P_W(\theta)|\Big|\sum^k_{j=0} e(j\ell\theta)\Big|d\theta\leq C\log (2+k).\Vert P_W\Vert_1.\tag 2.13
$$
Hence, from (2.5)
$$
\Vert  P_W\Vert_1 \lesssim \log (2+k_1) \Vert P_{W_1}\Vert_1 +\cdots+ \log (2+k_r)
 \Vert P_{W_r}\Vert_1.\tag 2.14
$$
Iteration of (2.14), using the geometric/arithmetic mean inequality and $n=o(\log|W|)$ by (2.2), (2.3), we get

\proclaim
{Lemma 4} The words $W\in\Cal W_n$ in the system (2.1) satisfy
$$
\Vert P_W\Vert_1< \Big(C\frac {\log|W|}n\Big)^n \ll |W|^\ve.\tag 2.15
$$
\endproclaim

\noindent
{\bf Remark.}

Returning to Lemma 2 and (1.61), we obtain a description of the word $B_n$ of the type (2.1) (with $n$ replaced by $r<c\log n)$, taking
$\Cal W_0=\{1^1, W_0\}$ and $\Cal W_s=\{W_s\}$.

Applying the bound (2.15) gives
$$
\Vert P_{B_n}\Vert_1<(\log h_n)^{c\log n}< n^{c\log n}.\tag 2.16
$$
On the other hand, the best unconditional bound (1.63)
$$
\Big\Vert \sum^{h_n}_{m=1} \mu(m) e(m\theta)\Big\Vert_\infty\ll_A h_n n^{-A} \text { (for all $A$)}
$$
together with (2.16) still falls short of implying a nontrivial estimate on (1.62).

We will develop a more refined analysis based on Vinogradov's bound that we recall next.

\proclaim
{Lemma 5} (Theorem 13.9 in \cite {I-K}).

Assume
$$
\Big|\theta-\frac aq\Big|\leq \frac 1{q^2} \text { with } (a, q)=1.\tag 2.17
$$
Then
$$
\Big|\sum_{m\leq x} \mu(m) e(m\theta)\Big|< c\Big(q^{\frac 12} x^{-\frac 12}+ q^{-\frac 12} +x^{-\frac 15}\Big)^{\frac 12} (\log x)^4x.\tag 2.18
$$
\endproclaim

There is the following consequence

\proclaim
{Lemma 6} Let $0<\tau<\frac 13$ and assume
$$
\Big|\theta -\frac aq\Big|\leq\frac 1{qx^{1-\tau}} \text { with } q\leq x^{1-\tau}, (a, q) =1.\tag 2.19
$$
Let $\beta =\theta-\frac aq$. Then
$$
\Big|\sum\mu(m) e(m\theta)\Big| <C\big((q+x |\beta|)^{-\frac 14}+ x^{-\frac\tau 4}\big)x (\log x)^4.\tag 2.20
$$
\endproclaim

\noindent
{\bf Proof.}

From (2.18)
$$
\Big|\sum^x_1\mu(m) e(m\theta)\Big|< C\big(x^{-\frac\tau 4}+q^{-\frac 14}\big) x (\log x)^4.\tag 2.21
$$
Assume $\beta\not=0$.

Next, take $a_1, q_1 \leq M=[\frac 2{|\beta|}]$ such that $(a_1, q_1)=1$ and $\Big|\theta-\frac {a_1}{q_1}\Big|\leq \frac 1{q_1 M}$.
Thus $|\beta|\leq \Big|\frac aq-\frac {a_1}{q_1}\Big|= \frac 1{q_1 M}\leq |\frac aq -\frac {a_1}{q_1}|+\frac \beta 2$ implying
$0\not= |\frac aq-\frac {a_1}{q_1}|\geq \frac 1{qq_1}$.

Hence $\frac 1{qq_1} \leq |\beta|+ \frac 1{q_1 M} < 2|\beta|$ and $q_1\geq \frac 1{2q|\beta|}$.

Apply (2.18) with $\frac aq$ replaced by $\frac {a_1}{q_1}$.
We obtain a second bound
$$
\align
\Big|\sum^x_1\mu(m)e(m\theta)\big|&< C \big(q_1^{\frac 14} x^{-\frac 14}+ q_1^{-\frac 14}+ x^{-\frac 1{10}}\big)(\log x)^4 x\\
&<C \big(|\beta|^{-\frac 14}x^{-\frac 14} +|\beta|^{\frac 14} q^{\frac 14}+x^{-\frac 1{10}}\big)(\log x)^4 x\\
&\overset{(2.19)}\to< C\big(|\beta|^{-\frac 14} x^{-\frac 14}+x^{-\frac 1{10}}\big)(\log x)^4x\tag 2.22
\endalign
$$
and (2.20) follows from (2.21), (2.22).

\bigskip

Let $W\in \Cal W_n, |W|=N$.
\bigskip

Define for $Q<N^{1-\tau}, K<N^\tau$ (dyadic) the sets
$$
V_Q =\bigcup_{\Sb (a, q)=1\\ q\sim Q\endSb} \Big[ \theta\in\Bbb T; \Big|\theta-\frac aq\Big|\leq \frac 1N\Big]\tag 2.23
$$
and
$$
V_{Q, K} =\bigcup_{\Sb (a, q)=1\\ q\sim Q\endSb} \Big[\theta\in\Bbb T; \Big|\theta-\frac aq\Big|\sim\frac KN\Big].\tag 2.24
$$
It follows from Lemma 6 that
$$
\Big\Vert \sum^N_{1} \mu(m) e(m\theta)\Big\Vert_{L^\infty (V_{Q, K)}} <C\big[(Q+K)^{-\frac 14}+N^{-\frac\tau 2}\big] N(\log N)^4.\tag 2.25
$$
First, from (2.15), (2.25)
$$
\sum_{\max(Q, K)>Q_0}\int_{V_{Q, K}} |P_W(\theta)| \, \Big|\sum^N_1 \mu (m)e(m\theta)\Big|d\theta \ll_\ve (Q_0^{-\frac 14}+N^{-\frac\tau 4})
N^{1+\ve}\tag 2.26
$$
and we therefore assume $Q, K< N^\ve$ in what follows.

Set $M=\frac NK$.  Note that if $V$ is a word of length $|V|=\ell<M$ and $\ell k>M$, we may write
$$
V^k=V_1^{k_1}V_2 \text {  where $V_1$ and $V_2$ are powers of $V$, $|V_1|\sim M$ and $|V_2|\leq M$}.\tag 2.27
$$
Write according to (2.1)
$$
W=W_1^{k_1}\cdots W_r^{k_r} \text { with } W_i\in\Cal W_{n-1}'.
$$
Fix $i =1, \ldots, r$ and distinguish the following cases.

\itemitem{(2.28)} $|W_i|>M$. Re-express then $W_i$ in words from $\Cal W_{n-2}'$.

\itemitem{(2.29)} $\frac M{k_i}<|W_i|\leq M$.

Write according to (2.27) $W_i^{k_i}=W_{i, 1}^k W_{i, 2}$ with $W_{i, 1}, W_{i, 2}$ powers of $W_i$ such that $|W_{i, 1}|\sim M$, $|W_{i, 2}|<M$.

We do not express $W_{i, 1}, W_{i, 2}$ in lower order words.

\itemitem{(2.30)} $|W_i|\leq\frac M{k_i}$.

Do not re-express $W_i$ in lower order words.

From triangle inequality and Lemma 3, estimate on $I_a =\big[|\theta-\frac aq|\sim\frac KN\big]$
$$
\align
\int_{I_a} |P_W|&\leq C\sum_{|W_i|>M} \log (k_i+2)\int_{I_a'}|P_{W_i}|\tag 2.31\\
&+C\sum_{\frac M{k_i} <|W_i|\leq M} \log\frac {k_i|W_i|}M \int_{I_a'} |P_{W_{i, 1}}|\tag 2.32
\\
&+ \sum_{\frac M{k_i} <|W_i|\leq M}\int_{I_a} |P_{W_{i, 2}}|\tag 2.33\\
&+\sum_{k_i |W_i|\leq M}\int_{I_a} |P_{{W_i}^{k_i}}|\tag 2.34\\
&+\sum_{|W_i|>M} \log (2+k_i) e^{-c(\log|W_i|)^{4/3}}\tag 2.35\\
&+\sum_{\frac M{k_i} <|W_i|\leq M} \log\Big(\frac {k_i|W_i|}M+2\Big) e^{-c(\log M)^{4/3}}\tag 2.36
\endalign
$$
where $I_a'$ is centered at $\frac aq$ of size $\frac {(\log N)^2}M$.

Since in (2.32), $\frac {k_i|W_i|}M\leq K$, (2.32), (2.33), (2.34) may be bounded by
$$
C\log K\int_{I_a'}|P_{W'}|\tag 2.37
$$
with $W'$ of the form $W'=V^k$, $V\in \Cal W_{n-1}'$ and $|W'|\lesssim M$.

Since $M>N^{1-\ve}$, clearly
$$
(2.35), (2.36)\lesssim \log K.e^{-c(\log N)^{4/3}}< e^{-c(\log N)^{4/3}}.\tag 2.38
$$
Repeat the preceding with each of the terms of (2.31) and iterate.

This leads to contributions of the form (with various $s$)
$$
\align
&\log (k_1+2)\cdots\log (k_s+2)\int_{I_a''}|P_{W'}|\tag 2.39\\
^+\\
&\log (k_1+2)\cdots \log (k_s+2) e^{-c(\log N)^{4/3}}\tag 2.40
\endalign
$$
where $W'=V^k, V\in\Cal W_{n-s}'$ and $|W'|\lesssim M$,
$$
k_1\ldots k_s \lesssim K\tag 2.41
$$
and $I_a''$ is centered at $\frac aq$ of size $s.\frac {(\log N)^2}M$.

By (2.2), (2.3), $s=o(\log K)$ and therefore
$$
\log(2+k_1) \ldots \log (2+ k_s)< \Big(C\frac {\log K}s\Big)^s <K^\ve.\tag 2.42
$$
From the preceding
$$
\sum_{\Sb (a, q)=1\\ q\sim Q\endSb} \int_{I_a} |P_W|< K^\ve\sum_{\Sb (a, q) =1\\ q\sim Q\endSb} 
\int_{I_a''}|P_{W'}|+ e^{-c(\log N)^{4/3}}\tag 2.43
$$
where $W'$ is of the form $W'=V^k$, $V\in \Cal W_{n-1}'$, $|W'|\lesssim M$ and $I_a''$ centered at $\frac aq$ of size $|I_a''|=\frac {(\log N)^3}M$.
\medskip

Assume $|W'|<Q^2$. Estimate
$$
\sum_{\Sb (a, q)=1\\ q\sim Q\endSb} \int_{I_a''} |P_{W'} | \leq \frac {(\log N)^3}M
\sum_{\Sb (a, q)=1\\ q\sim Q\endSb} \Big|P_{W'} \Big(\theta+\frac aq\Big)\Big|
\quad \text { for some $\theta$}.\tag 2.44
$$
Let $\ell=|W'|$ and take $K$ satisfying (2.8)-(2.10).
Hence $P_{W'}=P_{W'}*K$ and
$$
\sum_{\Sb (a, q)=1\\ q\sim Q\endSb} \Big|P_{W'} \Big(\theta+\frac aq\Big)\Big| \leq \int|P_{W'}(\eta)|\Big[\sum_{\Sb (a, q)=1\\ q\sim Q\endSb}
\Big|K\Big(\theta-\eta+\frac aq\Big)\Big|\Big] d\eta.
$$
From (2.10) and the separation $\Vert\frac aq-\frac {a'}{q'}\Vert \gtrsim \frac 1{Q^2}$ if $\frac aq\not= \frac {a'} {q'}, q, q'\sim Q$, it
follows that
$$
\align
\sum_{\Sb (a, q) =1\\ q\sim Q\endSb} \Big|K\Big(\theta-\eta +\frac aq\Big)\Big| &\lesssim \ell\sum_{\Sb (a, q)=1\\ q\sim Q\endSb}
e^{-c(\ell\Vert\theta -\eta+\frac aq\Vert)^{2/3}}\\
&\lesssim \ell \sum_{j=0}^{Q^2} e^{-c(\frac {\ell j}{Q^2})^{2/3}}\\
&\lesssim Q^2.
\endalign
$$
Thus
$$
(2.44)\lesssim \frac {(\log N)^3}M Q^2 \Vert P_{W'}\Vert_1
$$
and applying Lemma 4 to $W'\in\Cal W_{n'}$
$$
\lesssim \frac {(\log N)^3}M Q^2 |W'|^\ve < \frac {(\log N)^3}M Q^{2+\ve}.\tag 2.45
$$
If $|W'|>Q^2$, decompose in lower order words using (2.1) until obtaining words of size at most $Q^2$.
From triangle inequality and (2.45)
$$
\align
\sum_{\Sb (a, q)=1\\ q\sim Q\endSb} \int_{I_a''}| P_{W'}|& \leq \sum_\alpha\sum_{\Sb (a, q)=1\\ q\sim Q\endSb} 
\int_{I_a''}| P_{W_\alpha}| \text { with } |W| =\sum |W_\alpha|, |W_\alpha| \lesssim Q^2\\
&<\frac {(\log N)^3}M Q^2 \sum_\alpha |W_\alpha|^\ve.\tag 2.46
\endalign
$$

Note that since each word $W_\alpha$ occurs at previous stage in a word of size at least $Q^2$, 
the number of $\alpha$'s in (2.46) is at most
$C\frac {|W'|}{Q^2}\lesssim \frac M{Q^2}$.  Therefore
$$
(2.36)\leq \frac {(\log N)^3}M Q^2 \frac M{Q^2} Q^\ve<(\log N)^3 Q^\ve.\tag 2.47
$$
From (2.43), (2.47),
$$
\int_{V_{Q, K}} |P_W|\lesssim (\log N)^3 Q^\ve K^\ve.\tag 2.48
$$
Combined with (2.25), we proved that for $Q, K<N^\ve$
$$
\align
\int_{V_{Q, K}} |P_W(\theta)| \, \Big|\sum_1^N \mu(n)e(n\theta)\Big|d\theta&< c(\log N)^7 Q^\ve K^\ve(Q+K)^{-\frac 14}N\\
& <c(\log N)^7 (Q+K)^{-\frac 15}N.\tag 2.49
\endalign
$$
Since $\Vert\sum^N_1\mu(n)e(n\theta)\Vert_\infty \ll_A(\log N)^{-A}N$, also
$$
\int_{V_{Q, K}} |P_Q(\theta)|\Big|\sum^N_1\mu(n)e(n\theta)\Big|d\theta\ll N.\frac {Q^2K}N(\log N)^{-A}N<Q^2K(\log N)^{-A}N.\tag 2.50
$$
From (2.50) and summation of (2.49) over dyadic ranges of $Q, K$, we obtain

\proclaim
{Theorem 7} Let $\{\Cal W_n; n\geq 1\}$ be a symbolic system with properties (2.1)-(2.3) as described in the beginning of this section.
Then, if $W\in \bigcup \Cal W_n$ and $|W|=N$, we have
$$
\int_{\Bbb T} |P_W(\theta)|\, \Big|\sum^N_1 \mu(m)e(m\theta)\Big|d\theta \ll_A N(\log N)^{-A}.\tag 2.51
$$
\endproclaim

Recalling the Remark after Lemma 4, Theorem 2 completes in particular the proof of Theorem 1.

\bigskip

\noindent
{\bf (3). Further Comments and an Application to Interval Exchange Transformations}

\noindent
{\bf (1).} Theorem 1 and related discussions apply equally well to the Liouville function.

\noindent
{\bf (2).} Theorem 2 has other applications.
In particular, it allows us to prove the Moebius disjointness property for the following `classical' rank-one systems (cf. \cite{Fe}).

\itemitem {(i).} Generalized Chacon systems defined symbolically by
$$
B_{n+1} =B_n^{p_n} 1 B_n^{q_n}\tag 3.1
$$
with $p_n+q_n\to\infty$.

\itemitem {(ii).} Katok's systems
$$
B_{n+1} =B_n^{p_n}(B_n 1)^{p_n}\tag 3.2
$$
with $p_n\to\infty$ fast enough.

Both examples are rigid and weakly mixing. Katok's map appears
as a special case of a three-interval exchange transformation (cf. \cite{F-H-Z$_4$}) that will be discussed later
in greater detail in this section.

It is indeed easily verified that the above systems satisfy the condition for Theorem 2 to apply

\noindent
{\bf (3)}. For later discussion, it will be useful to remove the logarithmic factor in (2.48).
This may be achieved in the following way.

First, (2.6) may be stated in the form $(|W|=\ell)$
$$
\int_I| P_W(\theta)| \Big|\sum^k_{j=0} e(j\theta)\Big| d\theta \leq C\log (2+k)\int_{I'} | P_W|+C\ell|I| \log (2+k)e^{-cB^{2/3}}\tag 3.3
$$
with $I'=I+[-\frac B\ell, \frac B\ell]$ and $B>1$ a parameter.

This follows easily by an inspection of the proof of Lemma 3.

Take $B=(\log K)^2$.
The expression (2.35) is replaced by
$$
\sum_{|W_i|>M} \log (2+k_i) e^{-c(\log K)^{4/3}} \frac {|W_i|}M\lesssim \log (2+k_i) e^{-c(\log K)^{4/3}}
$$
and similarly for (2.36); (2.40) becomes
$$
\log (k_1 +2)\cdots \log (k_s+2) e^{-c(\log K)^{4/3}} < K e^{-c(\log K)^{4/3}} < e^{-c(\log K)^{4/3}}$$
while in (2.43), $|I_a''|<\frac {(\log K)^3}M$ with an 
error term bounded by $e^{-c(\log K)^{4/3}}$.

Following the rest of the argument verbatim, this leads eventually to an estimate
$$
\int_{V_{Q, K}}| P_W|\lesssim (\log K)^3 Q^\ve K^\ve \ll Q^\ve K^\ve\tag 3.4
$$
instead of (2.48).
\medskip

\noindent
{\bf (4).} The interest of (3.4) is that it enables us to exploit the usual Hardy-Littlewood circle method to study sums over the primes.
Assume \hfill\break $W=(x_1, \ldots, x_N)\in \{0, 1\}^N$ satisfies (3.4).
Fix $q\in\Bbb Z_+$ a large parameter (independent of $N$).
Then one can show in particular that
$$
\sum^N_1 x_j\Lambda(j+n) =\frac q{\phi(q)}\sum^N_1 x_j1_{\pi_q (j+n)\in (\Bbb Z/q\Bbb Z)^*} +O\Big(\frac {N+n}{(\log q)^{1/2}}\Big).\tag 3.5
$$
If furthermore $W$ is a word produced by a weakly mixing and uniquely ergodic system, the first term in (3.5) equals (since $q$ is fixed)
$$
\sum^N_1 x_j+o(N).
$$
Indeed, one has for any $a\in \Bbb Z/q\Bbb Z$
$$
\align
\sum^N_1 x_j 1_{\pi_q(j)=a}&= \frac 1q \sum^{q-1}_{k=0} \sum^N_{j=1}  e_q\big(k(j-a)\big) x_j\\
&=\frac 1q \sum^N_{j=1} x_j+(3.6)\\
\endalign
$$
where
$$
|(3.6)|\leq \max_{0< k< q} \Big|\sum^N_{j=1} e_q (kj)x_j\Big|\overset {N\to \infty}\to\longrightarrow 0.
$$
Hence
$$
\sum^N_1 x_j\Lambda(j+n)=\sum^N_1x_j+ O\Big(\frac {N+n}{(\log q)^{1/2}}\Big)+o(N).\tag 3.7
$$

\noindent
{\bf (5).}  Beyond the examples (3.1), (3.2), there is a natural family of systems that fit in the frame work discussed in \S2 and above,
$n\ell$ a rather large class of three-interval exchange transformations (3-IET).
We recall a few facts referring to the papers \cite{F-H-Z$_{1,2,3,4}$}.

Given $\alpha, \beta>0, \alpha+\beta<1$, define a transformation $T$ on $[0, 1]$  by
$$
Tx= \left\{ \aligned &x+1-\alpha\text { if } x\in [0, \alpha[\\
&x+1-2\alpha \text { if } x\in [\alpha, \alpha +\beta[\\
&x-\alpha- \beta \text { if } x\in [\alpha +\beta, 1[.
\endaligned \right.
\tag 3.8
$$
We assume $\alpha, \beta, 1$ independent over $\Bbb Q$, implying Keane's infinite distinct orbit condition, which in turn ensures minimality and unique
ergodicity of $T$.

Next, one associates to $(\alpha, \beta)$ a sequence $(n_k, m_k, \ve_{k+1})_{k\geq 0}$ with $n_k, m_k\in\Bbb Z_+$ and $\ve_{k+1}=\pm
1$ \big(the three-interval
expansion of $(\alpha,\beta)$\big).
The system $T$ may then be described symbolically using three return words $A_k, B_k, C_k$ satisfying the recursive relations
$$
\left\{
\aligned
& A_k= A_{k-1}^{n_k-1} C_{k-1} B_{k-1}^{m_k-1} A_{k-1}\\
&B_k= A_{k-1}^{n_k-1} C_{k-1} B_{k-1}^{m_k}\\
&C_k = A_{k-1}^{n_k-1} C_{k-1} B_{k-1}^{m_k-1} 
\endaligned\right.\tag 3.9
$$
if $\ve_{k+1}=1$ and
$$
\left\{
\aligned
&A_k= A_{k-1}^{n_k-1} C_{k-1} B_{k-1}^{m_k}\\
&B_k =A_{k-1}^{n_k-1} C_{k-1} B_{k-1}^{m_k-1} A_{k-1}\\
&C_k = A_{k-1}^{n_h-1} C_{k-1} B_{k-1}^{m_k} A_{k-1}
\endaligned
\right.
\tag 3.10
$$
if $\ve_{k+1} =-1$, 

\noindent
with initial words $A_0, B_0, C_0$ satisfying $\big||A_0|-|B_0|\big|=1$.

Let $a_k=|A_k|, b_k= |B_k|, c_k=|C_k|$.
Note that $|a_k-b_k|= |a_{k-1} -b_{k-1}|=1$ and $c_k \leq 2a_k$.

According to \cite {F-H-Z$_3$} (Theorem 3.5), a sufficient condition to ensure weak-mixing is that
$$
\int_{k} \frac {\min(m_k, n_k)}{m_k+n_k}>0.\tag 3.11
$$
In order to fulfill moreover the main assumption from \S2 \big(see also the Remark following (2.2), (2.3)\big), assume also that
$$
\min (m_k, n_k)> C_0\text { for all $k$}\tag 3.12
$$
where $C_0$ is a sufficiently large constant.

Then (3.7) will hold for words $W(x_1, \ldots, x_N)$ in
$$
\bigcup_k\{ A_k^m, B_k^m, C_k^m; m\in\Bbb Z_+\}.\tag 3.15
$$
Given an arbitrary word $W=(x_1, \ldots, x_N)$ in the language of the system, one may approximate $W$ by a collection of shifts of words $W'$ in  (3.13) 
of size $|W'|\sim |W|$ (depending on the approximation).
Applying $W'$ to each shifted word $W'$, it follows that
$$
\sum^N_1 \Lambda(j) x_j =\sum^N_1 x_j+o(N),\tag 3.14
$$
Hence, we proved

\proclaim
{Theorem 8}
Assume that $T$ is a 3-IET satisfying the Keane condition, (3.11) and (3.12).
Then $T$ satisfies the Moebius disjointness property and also a prime number theorem.
\endproclaim

\enddocument